\documentclass[graybox]{svmult}
\usepackage{graphicx}
\usepackage{enumerate}
\usepackage{amsfonts}
\usepackage{pstricks}
\usepackage{verbatim}
\usepackage{bbm}
\usepackage{bbold}
\usepackage{mathrsfs, hyphenat}
\def\int{\mathbb Z}

\def\ZN{{\mathbb Z}_N}
\def\H{{\mathbb H}}
\def\Z2{{\mathbb Z}_2}

\def\E_B{{\frak E}_{{\frak B}}}%
\def\Lin_F{\mbox{Lin}(\frak F)}
\def\1{\mathbbm 1}%
\def\module_1{{\frak C}_N^{\bar 1}}%
\def\A{{\mathscr{A}}}%
\usepackage{mathptmx}       
\usepackage{helvet}         
\usepackage{courier}        
\usepackage{type1cm}        
\usepackage{makeidx}         
\usepackage{graphicx}        
\usepackage{multicol}        
\usepackage[bottom]{footmisc}

\usepackage{amssymb}
\usepackage{amsfonts}
\usepackage{comment}
\usepackage{mathtools}

\makeindex             

\newcommand{\dfn}[1]{\begin{definition} #1 \end{definition}}
\newcommand{\thm}[1]{\begin{theorem} #1 \end{theorem}}
\newcommand{\prop}[1]{\begin{proposition} #1 \end{proposition}}

\def\reducedqp{\mathscr A_{q}[x,y]}
\def\reducedqpk{{\mathscr A}^{k}_{q}[x,y]}

\def\reducedqpzero{{\mathscr A}^{0}_{q}[x,y]}
\def\reducedqpone{{\mathscr A}^{1}_{q}[x,y]}
\def\reducedqpk{{\mathscr A}^{k}_{q}[x,y]}
\def\reducedqpkplusone{{\mathscr A}^{k+1}_{q}[x,y]}

\begin{document}

\title*{Semi-Commutative Galois Extension and Reduced Quantum Plane}
\author{Viktor Abramov, Md. Raknuzzaman}
\institute{%
  Viktor Abramov \at
  Institute of Mathematics, University of Tartu, \\
  Liivi 2--602, Tartu 50409, Estonia, \\
  \email{viktor.abramov@ut.ee}
  \and
  Md. Raknuzzaman\at
  Institute of Mathematics, University of Tartu, \\
  Liivi 2--602, Tartu 50409, Estonia, \\
  \email{raknuzza@ut.ee}
}
%
%
\maketitle

\abstract*{In this paper we show that a semi-commutative Galois extension of associative unital algebra by means of an element $\tau$, which satisfies $\tau^N=\mathbb 1$ ($\mathbb 1$ is the identity element of an algebra and $N\geq2$ is an integer) induces a structure of graded $q$-differential algebra, where $q$ is a primitive $N$th root of unity. A graded $q$-differential algebra with differential $d$, which satisfies $d^N=0, N\geq2$, can be viewed as a generalization of graded differential algebra. The subalgebra of elements of degree zero and the subspace of elements of degree one of a graded $q$-differential algebra together with a differential $d$ can be considered as a first order noncommutative differential calculus. In this paper we assume that we are given a semi-commutative Galois extension of associative unital algebra, then we show how one can construct the graded $q$-differential algebra and when this algebra is constructed we study its first order noncommutative differential calculus. We also study the subspaces of graded $q$-differential algebra of degree greater than one which we call the higher order noncommutative differential calculus induced by a semi-commutative Galois extension of associative unital algebra. Finally we show that a reduced quantum plane can be viewed as a semi-commutative Galois extension of a fractional one-dimensional space and we apply the noncommutative differential calculus developed in the previous sections to a reduced quantum plane. }

\abstract{In this paper we show that a semi-commutative Galois extension of associative unital algebra by means of an element $\tau$, which satisfies $\tau^N=\mathbb 1$ ($\mathbb 1$ is the identity element of an algebra and $N\geq2$ is an integer) induces a structure of graded $q$-differential algebra, where $q$ is a primitive $N$th root of unity. A graded $q$-differential algebra with differential $d$, which satisfies $d^N=0, N\geq2$, can be viewed as a generalization of graded differential algebra. The subalgebra of elements of degree zero and the subspace of elements of degree one of a graded $q$-differential algebra together with a differential $d$ can be considered as a first order noncommutative differential calculus. In this paper we assume that we are given a semi-commutative Galois extension of associative unital algebra, then we show how one can construct the graded $q$-differential algebra and when this algebra is constructed we study its first order noncommutative differential calculus. We also study the subspaces of graded $q$-differential algebra of degree greater than one which we call the higher order noncommutative differential calculus induced by a semi-commutative Galois extension of associative unital algebra. We also study the subspaces of graded $q$-differential algebra of degree greater than one which we call the higher order noncommutative differential calculus induced by a semi-commutative Galois extension of associative unital algebra. Finally we show that a reduced quantum plane can be viewed as a semi-commutative Galois extension of a fractional one-dimensional space and we apply the noncommutative differential calculus developed in the previous sections to a reduced quantum plane. }

\section{Introduction}
Let us briefly remind a definition of noncommutative Galois extension \cite{Kerner-Suzuki,Lawrynowicz-al-1,Lawrynowicz-al-2,Trovon}. Suppose $\tilde\A$ is an associative unital $\mathbb C$-algebra, $\A\subset \tilde\A$ is its subalgebra, and there is an element $\tau\in\tilde\A$ which satisfies $\tau \notin \A,\;\tau^N=\1$, where $N\geq 2$ is an integer and $\1$ is the identity element of $\tilde\A$. A noncommutative Galois extension of $\A$ by means of $\tau$ is the smallest subalgebra $\A[\tau]\subset \tilde\A$ such that $\A\subset \A[\tau]$, and $\tau\in\A[\tau]$. It should be pointed out that a concept of noncommutative Galois extension can be applied not only to associative unital algebra with a binary multiplication law but as well as to the algebra with a ternary multiplication law, for instant to a ternary analog of Grassmann and Clifford algebra \cite{Abramov-0,Lawrynowicz-al-2,Trovon}, and this approach can be used in particle physics to construct an elegant algebraic model for quarks.

A graded $q$-differential algebra can be viewed as a generalization of a notion of graded differential algebra if we use a more general equation $d^N=0, N\geq 2$ than the basic equation $d^2=0$ of a graded differential algebra. This idea was proposed and developed within the framework of noncommutative geometry \cite{Kapranov}, where the author introduced the notions of $N$-complex, generalized cohomologies of $N$-complex and making use of an $N$th primitive root of unity constructed an analog of an algebra of differential forms in $n$-dimensional space with exterior differential satisfying the relation $d^N=0$. Later this idea was developed in the paper \cite{Dubois-Violette-Kerner}, where the authors introduced and studied a notion of graded $q$-differential algebra. It was shown \cite{Abramov-1,Abramov-2,Abramov-Liivapuu,Abramov-Liivapuu-2} that a notion of graded $q$-differential algebra can be applied in noncommutative geometry in order to construct a noncommutative generalization of differential forms and a concept of connection.

In this paper we will study a special case of noncommutative Galois extension which is called a semi-commutative Galois extension. A noncommutative Galois extension is referred to as a semi-commutative Galois extension \cite{Trovon} if for any element $x\in\A$ there exists an element $x^\prime\in\A$ such that $x\,\tau=\tau\,x^\prime$. In this paper we show that a semi-commutative Galois extension can be endowed with a structure of a graded algebra if we assign degree zero to elements of subalgebra $\A$ and degree one to $\tau$. This is the first step on a way to construct the graded $q$-differential algebra if we are given a semi-commutative Galois extension. The second step is the theorem which states that if there exists an element $v$ of graded associative unital $\mathbb C$-algebra which satisfies the relation $v^N=\mathbb 1$ then this algebra can be endowed with the structure of graded $q$-differential algebra. We can apply this theorem to a semi-commutative Galois extension because we have an element $\tau$ with the property $\tau^N=\mathbb 1$, and this allows us to equip a semi-commutative Galois extension with the structure of graded $q$-differential algebra. Then we study the first and higher order noncommutative differential calculus induced by the $N$-differential of graded $q$-differential algebra. We introduce a derivative and differential with the help of first order noncommutative differential calculus developed in the papers \cite{Abramov-Kerner,Borowiec-Kharchenko}. We also study the higher order noncommutative differential calculus and in this case we consider a differential $d$ as an analog of exterior differential and the elements of higher order differential calculus as analogs of differential forms. Finally we apply our calculus to reduced quantum plane \cite{Coquereaux1}.





\section{Graded $q$-Differential Algebra Structure of Noncommutative Galois Extension}
In this section we remind a definition of noncommutative Galois extension, semi-commutative Galois extension, and show that given a semi-commutative Galois extension we can construct the graded $q$-differential algebra.

\vskip.3cm
\noindent
First of all we remind a notion of a noncommutative Galois extension \cite{Kerner-Suzuki,Lawrynowicz-al-1,Lawrynowicz-al-2,Trovon}.
\dfn{
Let $\tilde\A$ be an associative unital $\mathbb C$-algebra and $\A\subset \tilde\A$ be its subalgebra. If there exist an element $\tau\in\tilde\A$ and an integer $N\geq 2$ such that
\begin{itemize}
\item[i)] $\tau^N=\pm\1$,
\item[ii)] $\tau^k\notin\A$ for any integer $1\leq k\leq N-1$,
\end{itemize}
then the smallest subalgebra $\A[\tau]$ of $\tilde\A$ which satisfies
\begin{itemize}
\item[iii)] $\A\subset \A[\tau]$,
\item[iv)] $\tau\in \A[\tau]$,
\end{itemize}
is called the noncommutative Galois extension of $\A$ by means of $\tau$.
}
\noindent
In this paper we will study a particular case of a noncommutative Galois extension which is called a semi-commutative Galois extension \cite{Trovon}. A noncommutative Galois extension is referred to as a semi-commutative Galois extension if for any element $x\in\A$ there exists an element $x^\prime\in\A$ such that $x\,\tau=\tau\,x^\prime$. We will give this definition in terms of left and right $\A$-modules generated by $\tau$. Let $\A^1_{\mbox{l}}[\tau]$ and $\A^1_{\mbox{r}}[\tau]$ be respectively the left and right $\A$-modules generated by $\tau$. Obviously we have
$$
\A^1_{\mbox{l}}[\tau]\subset \A[\tau],\;\;\A^1_{\mbox{r}}[\tau]\subset\A[\tau].
$$
\dfn{
A noncommutative Galois extension $\A[\tau]$ is said to be a right (left) semi-commutative Galois extension if $\A^1_{\mbox{r}}[\tau]\subset \A^1_{\mbox{l}}[\tau]$ ($\A^1_{\mbox{l}}[\tau]\subset \A^1_{\mbox{r}}[\tau]$). If $\A^1_{\mbox{r}}[\tau]\equiv \A^1_{\mbox{l}}[\tau]$ then a noncommutative Galois extension will be referred to as a semi-commutative Galois extension, and in this case $\A^1[\tau]=\A^1_{\mbox{r}}[\tau]=\A^1_{\mbox{l}}[\tau]$ is the $\A$-bimodule.
}
\noindent
It is well known that a bimodule over an associative unital algebra $\A$ freely generated by elements of its basis induces the endomorphism from an algebra $\A$ to the algebra of square matrices over $\A$. In the case of semi-commutative Galois extension we have only one generator $\tau$ and it induces the endomorphism of an algebra $\A$. Indeed let $\A[\tau]$ be a semi-commutative Galois extension and $\A^1[\tau]$ be its $\A$-bimodule generated by $[\tau]$. Any element of the right $\A$-module $\A^1_{\mbox{r}}[\tau]$ can be written as $\tau\,x$, where $x\in \A$. On the other hand $\A[\tau]$ is a semi-commutative Galois extension which means $\A^1_{\mbox{r}}[\tau]\equiv \A^1_{\mbox{l}}[\tau]$, and hence each element $x\,\tau$ of the left $\A$-module can be expressed as $\tau\,\phi_\tau(x)$, where $\phi_\tau(x)\in\A$. It is easy to verify that the linear mapping $\phi:x\to \phi_\tau(x)$ is the endomorphism of subalgebra $\A$, i.e. for any elements $x,y\in \frak A$ we have $\phi_\tau(xy)=\phi_\tau(x)\phi_\tau(y)$. This endomorphism will play an important role in our differential calculus, and in what follows we will also use the notation $\phi_\tau(x)=x_\tau$. Thus
$$
u\,\tau=\tau\,\phi_\tau(x), \quad u\,\tau=\tau\,u_\tau.
$$
It is clear that
$$
\phi_\tau^N=\mbox{id}_{\A},\,\,\,u_{\tau^N}=u,
$$
because for any $u\in\A$ it holds $u\,\tau^N=\tau^N\,\phi^N(u)$ and taking into account that $\tau^N=\1$ we get $\phi_\tau^N(u)=u$.
\prop{
Let $\A[\tau]$ be a semi-commutative Galois extension of $\A$ by means of $\tau$, and $\A_{\mbox{l}}^k[\tau],\A_{\mbox{r}}^k[\tau]$ be respectively the left and right $\A$-modules generated by $\tau^k$, where $k=1,2,\ldots, N-1$. Then  $\A_{\mbox{l}}^k[\tau]\equiv \A_{\mbox{r}}^k[\tau]=\A^k[\tau]$ is the $\A$-bimodule, and
$$
\A[\tau]=\oplus_{k=0}^{N-1}\A^k[\tau]=\A^0[\tau]\oplus\A^1[\tau]\oplus\ldots\oplus\A^{N-1}[\tau],
$$
where $\A^0[\tau]\equiv \A$.
\label{proposition Galois extension}
}
\noindent
Evidently the endomorphism of $\A$ induced by the $\A$-bimodule structure of $A^k[\tau]$ is $\phi^k$, where $\phi:\A\to \A$ is the endomorphism induced by the $\A$-bimodule $\A^1[\tau]$. We will also use the notation $\phi^k(x)=x_{\tau^k}$.

\vskip.3cm
\noindent
It follows from Proposition \ref{proposition Galois extension} that a semi-commutative Galois extension $\A[\tau]$ has a natural $\ZN$-graded structure which can be defined as follows: we assign degree zero to each element of subalgebra $\A$, degree 1 to $\tau$ and extend this graded structure to a semi-commutative Galois extension $\A[\tau]$ by determining the degree of a product of two elements as the sum of degree of its factors. The degree of a homogeneous element of $\A[\tau]$ will be denoted by $|\;\;|$. Hence $|u|=0$ for any $u\in\A$ and $|\tau|=1$.

\vskip.3cm
\noindent
Now our aim is to show that given a noncommutative Galois extension we can construct a graded $q$-differential algebra, where $q$ is a primitive $N$th root of unity. First of all we remind some basic notions, structures and theorems of theory of graded $q$-differential algebras.

\vskip.3cm
\noindent
Let $\A=\oplus_{k\in{\mathbb Z}_N}{\A}^k=\A^0\oplus\A^1\oplus\ldots\oplus\A^{N-1}$ be a ${\mathbb Z}_N$-graded associative unital $\mathbb C$-algebra with identity element denoted by $ \1$. Obviously the subspace $\A^0$ of elements of degree 0 is the subalgebra of a graded algebra $\A$.
Every subspace $\A^k$ of homogeneous elements of degree $k\geq 0$ can be viewed as the $\A^0$-bimodule.
The graded $q$-commutator of two homogeneous elements $u,v\in \A$ is defined  by
$$
[v,u]_q=v\,u-q^{|v||u|}u\,v.
$$
A graded $q$-derivation of degree $m$ of a graded algebra $\A$ is a linear mapping $d:\A\to\A$ of degree $m$, i.e. $d:\A^i\to \A^{i+m}$, which satisfies the graded $q$-Leibniz rule
\begin{equation}
d(u\,v)=d(u)\,v+q^{ml}u\,d(v),
\label{graded q-Leibniz}
\end{equation}
where $u$ is a homogeneous element of degree $l$, i.e. $u\in \A^l$. A graded $q$-derivation $d$ of degree $m$ is called an inner graded $q$-derivation of degree $m$ induced by an element $v\in\A^m$ if
\begin{equation}
d(u)=[v,u]_q=v\,u-q^{ml}u\,v,
\end{equation}
where $u\in \A^l$.

\vskip.3cm
\noindent
Now let $q$ be a primitive $N$th root of unity, for instant $q=e^{2\pi i/N}$. Then
$$
q^N=1,\;\; 1+q+\ldots+q^{N-1}=0.
$$
A graded $q$-differential algebra
is a graded associative unital algebra $\A$ endowed with a graded $q$-derivation $d$ of degree one which satisfies $d^N=0$. In what follows a graded $q$-derivation $d$ of a graded $q$-differential algebra $\A$ will be referred to as a graded $N$-differential. Thus a graded $N$-differential $d$ of a graded $q$-differential algebra is a linear mapping of degree one which satisfies a graded $q$-Leibniz rule and $d^N=0$. It is useful to remind that a graded differential algebra is a graded associative unital algebra equipped with a differential $d$ which satisfies the graded Leibniz rule and $d^2=0$. Hence it is easy to see that a graded differential algebra is a particular case of a graded $q$-differential algebra when $N=2, q=-1$, and in this sense we can consider a graded $q$-differential algebra as a generalization of a concept of graded differential algebra. Given a graded associative algebra $\A$ we can consider the vector space of inner graded $q$-derivations of degree one of this algebra and put the question: under what conditions an inner graded $q$-derivation of degree one is a graded $N$-differential? The following theorem gives  answer to this question.

\thm{
Let $\A$ be a $\mathbb Z_N$-graded associative unital $\mathbb C$-algebra and $d(u)=[v,u]_q$ be its inner graded $q$-derivation induced by an element $v\in\A^1$. The inner graded $q$-derivation $d$ is the $N$-differential, i.e. it satisfies $d^N=0$, if and only if $v^N=\pm \1$.
}
\noindent
Now our goal is apply this theorem to a semi-commutative Galois extension to construct a graded $q$-differential algebra with $N$-differential satisfying $d^N=0$.
\prop{
Let $q$ be a primitive $N$th root of unity. A semi-commutative Galois extension $\A[\tau]$, equipped with the $\ZN$-graded structure described above and with the inner graded $q$-derivation $d\,=[\tau,\,]_q$ induced by $\tau$, is the graded $q$-differential algebra, and $d$ is its $N$-differential. For any element $\xi$ of semi-commutative Galois extension $\A[\tau]$ written as a sum of elements of right $\A$-modules $\A^k[\tau]$
$$
\xi=\sum_{k=0}^{N-1}\tau^k\,u_k=\1\,u_0+\tau\,u_1+\tau^2\,u_2+\ldots \tau^{N-1}\,u_{N-1},\;\;u_k\in \A,
$$
it holds
\begin{equation}
d \xi=\sum_{k=0}^{N-1}\tau^{k+1}(u_k-q^k\,(u_k)_\tau),
\label{semi-commutative graded q-differential}
\end{equation}
where $u_k\to (u_k)_\tau$ is the endomorphism of $\A$ induced by the bimodule structure of $\A^1[\tau].$
\label{proposition 2.5}
}
\section{First Order Differential Calculus over Associative Unital Algebra}

In this section we describe a first order differential calculus over associative unital algebra \cite{Borowiec-Kharchenko}. If an associative unital algebra is generated by a family of variables, which obey commutation relations, then one can construct a coordinate first order differential calculus over this algebra. A coordinate first differential calculus induces the partial derivatives with respect to generators of algebra and these partial derivatives satisfy the twisted Leibniz rule.

\vskip .3cm
\noindent A first order differential calculus is a triple $(\A, \mathscr M, d)$ where $\A$ is an associative unital algebra,  $\mathscr M$ is an $\A$-bimodule, and $d$, which is called a differential of first order differential calculus, is a linear mapping $d: \A \to \mathscr M$ satisfying the Leibniz rule $d(fh)=dfh+fdh$, where $f,h \in \A$.
A first order differential calculus $(\A, \mathscr M, d)$ is referred to as a coordinate first order differential calculus if an algebra $\A$ is generated by the variables $x^{1},x^{2},\cdots, x^{n}$ which satisfy the commutation relations, and an $\A$-bimodule $\mathscr M$, considered as a right $\A$-module, is freely generated by $dx^{1},dx^{2},\cdots, dx^{n}$. It is worth to mention that a first order differential calculus was developed within the framework of noncommutative geometry, and an algebra $\A$ is usually considered as the algebra of functions of a noncommutative space, the generators $x^{1},x^{2},\cdots, x^{n}$ of this algebra are usually interpreted as coordinates of this noncommutative space, and an $\A$-bimodule $\mathscr M$ plays the role of space of differential forms of degree one. In this paper we will use the corresponding terminology in order to stress a relation with noncommutative geometry.
\vskip .3cm


\noindent Let us consider a structure of coordinate first order differential calculus. This differential calculus induces the differentials $dx^{1},dx^{2},\cdots, dx^{n}$ of the generators $x^{1},x^{2},\cdots, x^{n}$. Evidently $dx^{1},dx^{2},\cdots, dx^{n}\in \mathscr M$. 
$\mathscr M$ is a bimodule, i.e. it has a structure of left $\mathscr A$-module and right $\mathscr A$-module. Hence for any two elements $f,h\in\A$ and $\omega\in \mathscr M$ it holds $(f\omega)h=f(\omega h)$. According to the definition of a coordinate first order differential calculus the right $\A$-module $\mathscr M$ is freely generated by the differentials of generators $dx^{1},dx^{2},\cdots,dx^{n}$. Thus for any $\omega \in \mathscr M$ we have $\omega=dx^{1}f_{1}+dx^{2}f_{2}+\cdots+dx^{n}f_{n}$ where $f_{1},f_{2},\cdots, f_{n}\in \A$. A coordinate first order differential calculus $(\A, \mathscr M, d)$ is an algebraic structure, which extends to noncommutative case the classical differential structure of a manifold. From the point of view of noncommutative geometry $\mathscr A$ can be viewed as an algebra of smooth functions, $d$ is the exterior differential, and $\mathscr M$ is the bimodule of differential $1$-forms. In order to stress this analogy we will call the elements of algebra $\mathscr A$ "functions" and the elements of $\A$-bimodule $\mathscr M$ ''$1$-forms".
\vskip .3cm
\noindent Because $\mathscr M$ is $\A$-bimodule, for any function $f\in \A$ we have two products $f\,dx^{i}$ and $dx^{i}\,f$. Since $dx^{1},dx^{2},\cdots,dx^{n}$ is the basis for the right $\mathscr A$-module $\mathscr M$, each element of $\mathscr M$ can be expressed as linear combination of $dx^{1},dx^{2},\cdots,dx^{n}$ multiplied by the functions from the right. Hence the element $fdx^{i}\in \mathscr M$ can be expressed in this way, i.e.
\begin{equation}
fdx^{i}=dx^{1} r^{i}_{1}(f)+dx^{2} r^{i}_{2}(f)+\ldots +dx^{n} r^{i}_{n}(f)=dx^{j} r^{i}_{j}(f),
\label{equation-1}
\end{equation}
where $r^{i}_{1}(f),r^{i}_{2}(f),\ldots,r^{i}_{n}(f)\in\A$ are the functions. Making use of these functions we can compose the square matrix
\begin{equation}
 R(f)=(r_{j}^{i}(f))=
\begin{pmatrix}
 r_{1}^{1}(f)& r_{1}^{2}(f)&\cdots& r_{1}^{n}(f)\\
\vdots&\vdots&\vdots&\vdots\\
 r_{n}^{1}(f)& r_{n}^{2}(f)&\cdots& r_{n}^{n}(f)
\end{pmatrix}.
\nonumber
\end{equation}
\noindent It is worth to point out that an entry $ r^{i}_{j}(f)$ stands on intersection of $i$-th column and $j$-th row.
This square matrix determines the mapping $R: \mathscr A \to  \text{Mat}_{n}(\A)$ where $ \text{Mat}_{n}(\A)$ is the algebra of $n$ order square matrices over an algebra $\A$. It can be proved
\vskip .3cm
\prop{
$R: \mathscr A \to  \text{Mat}_{n}(\A)$ is the homomorphism of algebras.
}
\begin{proof}
We need to prove that for any $f,g\in \A$ it holds $R(fg)=R(f)R(g)$.
Now according to the equation (\ref {equation-1}) we have
\begin{equation}
(fg)dx^{i}=dx^{j}r_{j}^{i}(fg).
\nonumber
\end{equation}
\vskip .3cm
\noindent The left hand side of the above relation can be written as
\begin{equation}
f(gdx^{i})=f(dx^{j}r_{j}^{i}(g))=(fdx^{j})r_{j}^{i}(g)=(dx^{k}r_{k}^{j}(f))r^{i}_{j}(g)=dx^{k}(r_{k}^{j}(f)r^{i}_{j}(g)).
\nonumber
\end{equation}
\noindent Now we can write
\begin{equation}
dx^{j}r_{j}^{i}(fg)=dx^{k}(r_{k}^{j}(f)r^{i}_{j}(g))\Rightarrow r_{k}^{i}(fg)=r_{k}^{j}(f)r^{i}_{j}(g),
\nonumber
\end{equation}
or in matrix form $R(fg)=R(f)R(g)$ which ends the proof.
\end{proof}
Let $\A,{\mathscr M},d$ be a coordinate first order differential calculus such that right $\A$-module $\mathscr M$ is a finite freely generated by the differentials of coordinates $\{dx_i\}_{i=1}^n$. The mappings $\partial_k:\A \to\A$, where $k\in\{1,2,\ldots,n\}$, uniquely defined by
\begin{equation}
df=dx^k\,\partial_k(f),\quad f\in \A,
\end{equation}
are called {the right partial derivatives} of a coordinate first order differential calculus. It can be proved
\begin{proposition}
If $\A,{\mathscr M},d$ is a coordinate first order differential calculus over an algebra $\A$ such that $\mathscr M$ is a finite freely generated right $\A$-module with a basis $\{dx_i\}_{i=1}^n$ then the right partial derivatives $\partial_k:\A \to\A$ of this differential calculus satisfy
\begin{equation}
\partial_k(fg)=\partial_k(f)\,g+r(f)^i_k\,\partial_i(g).
\label{twisted Leibniz rule for partial derivatives}
\end{equation}
\end{proposition}
\noindent
The property (\ref{twisted Leibniz rule for partial derivatives}) is called {the twisted (with homomorphism $R$) Leibniz rule}
\index{Twisted Leibniz rule for derivatives}
for partial derivatives.
\vskip.3cm
\noindent
If $\A$ is a graded $q$-differential algebra with differential $d$ then evidently the subspace of elements of degree zero $\A^0$ is the subalgebra of $\A$, the subspace of elements of degree one $\A^1$ is the $\A^0$-bimodule, a differential $d:\A^0\to \A^1$ satisfies the Leibniz rule. Consequently we have the first order differential calculus $(\A^0,\A^1,d)$ of a graded $q$-differential algebra $\A$. If $\A^0$ is generated by some set of variables then we can construct a coordinate first order differential calculus with corresponding right partial derivatives.

\section{First Order Differential Calculus of Semi-Commutative Galois Extension}
It is shown in Section 2 that given a semi-commutative Galois extension we can construct a graded $q$-differential algebra. In the previous section we described the structure of a coordinate first order differential calculus over an associative unital algebra, and at the end of this section we also mentioned that the subspaces $\A^0,\A^1$ of a graded $q$-differential algebra together with differential $d$ of this algebra can be viewed as a first order differential calculus over $\A^0$. In this section we apply an approach of first order differential calculus to a graded $q$-differential algebra of a semi-commutative Galois extension.

\vskip.3cm
\noindent
Let $\A[\tau]$ be a semi-commutative Galois extension of an algebra $\A$ by means of $\tau$. Thus we have an algebra $\A$ and $\A$-bimodule $\A^1[\tau]$. Next we have the $N$-differential $d:\A[\tau]\to\A[\tau]$ induced by $\tau$, and if we restrict this $N$-differential to the subalgebra $\A$ of Galois extension $\A[\tau]$ then $d:\A\to \A^1[\tau]$ satisfies the Leibniz rule. Consequently we have the first order differential calculus which can be written as the triple $(\A,d,\A^1[\tau])$. In order to describe the structure of this first order differential calculus we will need the vector space endomorphism $\Delta:\A\to\A$ defined by
$$
\Delta u=u-u_\tau,\quad u\in\A.
$$
For any elements $u,v\in\A$ this endomorphism satisfies
$$
\Delta(u\,v)=\Delta(u)\,v+u_\tau\,\Delta(v).
$$
Let us assume that there exists an element $x\in \A$ such that the element $\Delta x\in \A$ is invertible, and the inverse element will be denoted by $\Delta x^{-1}.$ The differential $dx$ of an element $x$ can be written in the form $dx=\tau\,\Delta x$ which clearly shows that $dx$ has degree one, i.e. $dx\in\A^1[\tau]$, and hence $dx$ can be used as generator for the right $\A$-module $\A^1[\tau]$. Let us denote by $\phi_{dx}:u\to \phi_{dx}(u)=u_{dx}$ the endomorphism of $\A$ induced by bimodule structure of $\A^1[\tau]$ in the basis $dx$. Then
\begin{equation}
u_{dx}=\Delta x^{-1}\,u_\tau\,\Delta x=\mbox{Ad}_{\Delta\,x}\,u_\tau.
\end{equation}
\dfn{
For any element $u\in \A$ we define the right derivative $\frac{du}{dx}\in \A$ (with respect to $x$)  by the formula
\begin{equation}
du=dx\,\frac{du}{dx}.
\label{derivative}
\end{equation}
}
\noindent
Analogously one can define the left derivative with respect to $x$ by means of the left $\A$-module structure of $\A^1[\tau]$. Further we will only use the right derivative which will be referred to as the derivative and often will be denoted by $u^\prime_x$. Thus we have the linear mapping
$$
\frac{d}{dx}:\A\to\A, \quad \frac{d}{dx}:u\mapsto u^\prime_x.
$$
\prop{
For any element $u\in \A$ we have
\begin{equation}
\frac{du}{dx}=\Delta x^{-1}\,\Delta u.
\end{equation}
The derivative (\ref{derivative}) satisfies the twisted Leibniz rule, i.e. for any two elements $u,v\in \A$ it holds
$$
\frac{d}{dx}(u\,v)=\frac{du}{dx}\,v+\phi_{dx}(u)\,\frac{dv}{dx}=\frac{du}{dx}\,v+\mbox{Ad}_{\Delta\,x}\,u_\tau\,\frac{dv}{dx}.
$$
\label{twisted Leibniz rule}
}
\noindent
We have constructed the first order differential calculus with one variable $x$, and it is natural to study a transformation rule of the derivative of this calculus if we choose another variable. From the point of view of differential geometry we will study a change of coordinate in one dimensional space. Let $y\in\A$ be an element of $\A$ such that $\Delta\,y=y-y_\tau$ is invertible.
\prop{
Let $x,y$ be elements of $\A$ such that $\Delta\,x,\Delta\,y$ are invertible elements of $\A$. Then
$$
dy=dx\,\,y^\prime_x,\;\; \frac{d}{dx}=y^\prime_x\,\frac{d}{dy},\;\;
     dx=dy\,\,x^\prime_y,\;\;\frac{d}{dy}=x^\prime_y\,\frac{d}{dx},
$$
where $x^\prime_y=(y^\prime_x)^{-1}$.
}
\noindent
Indeed we have $dy=\tau\,\Delta\,y,\,dx=\tau\,\Delta\,x$. Hence $\tau=dx\,\Delta\,x^{-1}$ and
$$
dy=dx\,(\Delta\,x^{-1}\Delta\,y)=dx\,y^\prime_x.
$$
If $u$ is any element of $\A$ the for the derivatives we have
$$
\frac{du}{dx}=\Delta\,x^{-1}\,\Delta\,u=(\Delta\,x^{-1}\,\Delta\,y)\,(\Delta\,y^{-1}\,\Delta\,u)=y^\prime_x\,\frac{du}{dy}.
$$
As an example of the structure of graded $q$-differential algebra induced by $d_\tau$ on a semi-commutative Galois extension we can consider the quaternion algebra $\mathbb H$. The quaternion algebra $\mathbb H$ is associative unital algebra generated over $\mathbb R$ by $i,j,k$ which are subjected to the relations
$$
i^2=j^2=k^2=-\1,\;i\,j=-j\,i=k,\;j\,k=-k\,j=i,\;k\,i=-i\,k=j,
$$
where $\1$ is the unity element of $\H$. Given a quaternion
$$
\frak q=a_0\,\1+a_1\,i+a_2\,j+a_3\,k
$$
we can write it in the form $\frak q=(a_0\,\1+a_2\,j)+i\,(a_1+a_3\,j)$. Hence if we consider the coefficients of the previous expression $z_0=a_0\,\1+a_2\,j,z_1=a_1+a_3\,j$ as complex numbers then $\frak q=z_0\,\1+i\,z_1$ which clearly shows that the quaternion algebra $\H$ can be viewed as the semi-commutative Galois extension $\mathbb C[i]$. Evidently in this case we have $N=2,q=-1$, and $\mathbb Z_2$-graded structure defined by $|\1|=0, |i|=1$. Hence we can use the terminology of superalgebras. It is easy to see that the subspace of odd elements (degree 1) can be considered as the bimodule over the subalgebra of even elements $a\,\1+b\,j$ and this bimodule induces the endomorphism $\phi:\mathbb C\to\mathbb C$, where $\phi(z)=\bar{z}$. Let $d$ be the differential of degree one (odd degree operator) induced by $i$. Then making use of \ref{semi-commutative graded q-differential} for any quaternion $\frak q$ we have
$$
d\frak q=d(z_0\,\1+i\,z_1)=-(\bar{z}_1+z_1)\,\1.
$$
Obviously $d^2\frak q=0.$
\section{Higher Order Differential Calculus of Semi-Commutative Galois Extension}
Our aim in this section is to develop a higher order differential calculus of a semi-commutative Galois extension $\A[\tau]$. This higher order differential calculus is induced by the graded $q$-differential algebra structure. In Section 2 it is mentioned that a graded $q$-differential algebra can be viewed as a generalization of a concept of graded differential algebra if we take $N=2,q=-1$. It is well known that one of the most important realizations of graded differential algebra is the algebra of differential forms on a smooth manifold. Hence we can consider the elements of the graded $q$-differential algebra constructed by means of a semi-commutative Galois extension $\A[\tau]$ and expressed in terms of differential $dx$ as noncommutative analogs of differential forms with exterior differential $d$ which satisfies $d^N=0$. In order to stress this analogy we will consider an element $x\in\A$ as analog of coordinate, the elements of degree zero as analogs of functions, elements of degree $k$ as analogs of $k$-forms, and we will use the corresponding terminology. It should be pointed out that because of the equation $d^N=0$ there are higher order differentials $dx,d^2x,\ldots,d^{N-1}x$ in this algebra of differential forms.

Before we describe the structure of higher order differentials forms it is useful to introduce the polynomials $P_k(x),Q_k(x)$, where $k=1,2,\ldots,N$. Let us remind that $\Delta x=x-x_\tau\in\A$. Applying the endomorphism $\tau$ we can generate the sequence of elements
$$
\Delta x_\tau=x_\tau-x_{\tau^2}, \Delta x_{\tau^2}=x_{\tau^2}-x_{\tau^3},\ldots, \Delta x_{\tau^{N-1}}=x_{\tau^{N-1}}-x.
$$
Obviously each element of this sequence is invertible. Now we define the sequence of polynomials $Q_1(x),Q_2(x),\ldots,Q_N(x)$, where
$$
Q_k(x)=\Delta x_{\tau^{k-1}}\Delta x_{\tau^{k-2}}\ldots \Delta x_{\tau}\Delta x.
$$
These polynomials can be defined by means of the recurrent relation
$$
Q_{k+1}(x)=(Q_k(x))_\tau\Delta x.
$$
It should be mentioned that $Q_k(x)$ is the invertible element and
$$
(Q_k(x))^{-1}=\Delta x^{-1}\Delta x^{-1}_\tau\ldots \Delta x^{-1}_{\tau^{k-1}}.
$$
We define the sequence of elements $P_1(x),P_2(x),\ldots,P_N(x)\in\A$ by the recurrent formula
$$
P_{k+1}(x)=P_{k}(x)-q^{k}\,(P_{k}(x))_\tau,\quad k=1,2,\ldots,N-1,
$$
and $P_1(x)=\Delta x$. Clearly $P_1(x)=Q_(x)$ and for the $k=2,3$ a straightforward calculation gives
\begin{eqnarray}
P_2(x) &=& x-(1+q)\,x_\tau+q\,x_{\tau^2},\nonumber\\
P_3(x) &=& x-(1+q+q^2)\,x_\tau+(q+q^2+q^3)\,x_{\tau^2}-q^3\,x_{\tau^3}.\nonumber
\end{eqnarray}
\prop{
If $q$ is a primitive $N$th root of unity then there are the identities
$$
P_{N-1}(x)+(P_{N-1}(x))_\tau+\ldots+(P_{N-1}(x))_{\tau^{N-1}}\equiv 0,\;\;\; P_N(x)\equiv 0.
$$
}
\noindent
Now we will describe the structure of higher order differential forms. It follows from the previous section that any 1-form $\omega$, i.e. an element of $\A^1[\tau]$, can be written in the form $\omega=dx\,u$, where $u\in\A$. Evidently $d:\A\to \A^1[\tau],\,d\omega=dx\,u^{\prime}_x$. The elements of $\A^2[\tau]$ will be referred to as 2-forms. In this case there are two choices for a basis for the right $\A$-module $\A^2[\tau]$. We can take either $\tau^2$ or $(dx)^2$ as a basis for $\A^2[\tau]$. Indeed we have
$$
(dx)^2=\tau^2\,Q_2(x).
$$
It is worth mentioning that the second order differential $d^2x$ can be used as the basis for $\A^2[\tau]$ only in the case when $P_2(x)$ is invertible. Indeed we have
$$
d^2x=\tau^2\,P_2(x),\quad d^2x=(dx)^2\,Q^{-1}_2(x)P_2(x).
$$
If we choose $(dx)^2$ as the basis for the module of 2-forms $\A^2[\tau]$ then any 2-form $\omega$ can be written as $\omega=(dx)^2\,u$, where $u\in\A$. Now the differential of any 1-form $\omega=dx\,u$, where $u\in\A$, can be expressed as follows
\begin{equation}
d\omega=(dx)^2\,\big( q\,u^\prime_x+Q_2^{-1}(x)P_2(x)\,u \big).
\label{differential of 1-form}
\end{equation}
It should be pointed out that the second factor of the right-hand side of the above formula resembles a covariant derivative in classical differential geometry. Hence we can introduce the linear operator $D:\A\to \A$ by the formula
\begin{equation}
Du=q\,u^\prime_x+Q_2^{-1}(x)P_2(x)\,u,\quad u\in\A.
\label{covariant derivative 1}
\end{equation}
If $\omega=dv, v\in \A$, i.e. $\omega$ is an exact form, then
$$
d\omega=d^2v=(dx)^2\,Dv^\prime_x=(dx)^2\,\big(q\,v^{\prime\prime}_x+Q_2^{-1}(x)P_2(x)\,v_x^\prime\big).
$$
If we consider the simplest case $N=2,q=-1$ then
$$
d^2v=0,\;\; P_2(x)\equiv 0,\;\; (dx)^2\neq 0,
$$
and from the above formula it follows that $v^{\prime\prime}_x=0$.
\prop{
Let $\A[\tau]$ be a semi-commutative Galois extension of algebra $\A$ by means of $\tau$, which satisfies $\tau^2=\1$, and $d$ be the differential of the graded differential algebra induced by an element $\tau$ as it is shown in Proposition \ref{proposition 2.5}. Let $x\in\A$ be an element such that $\Delta x$ is invertible. Then for any element $u\in\A$ it holds $u^{\prime\prime}_x=0$, where $u^\prime_x$ is the derivative (\ref{derivative}) induced by $d$. Hence any element of an algebra $\A$ is  linear with respect to $x$.
\label{linear}
}

\noindent
The quaternions considered as the noncommutative Galois extension of complex numbers (Section 3) provides a simple example for the above proposition. Indeed in this case $\tau=i, \A\equiv \mathbb C$, where the imaginary unit is identified with $j$, $(a\,\1+b\,j)_\tau=a\,\1-b\,j$. Hence we can choose $x=a\,\1+b\,j$ iff $b\neq 0$. Indeed in this case $\Delta x=x-x_\tau=a\,\1+b\,j-a\,\1+b\,j=2b\,j$, and $\Delta x$ is invertible iff $b\neq 0$. Now any $z=c\,\1+d\,j\in\A$ can be uniquely written in the form $z=\tilde c\,\1+\tilde d\,x$ iff
$$
\left|
  \begin{array}{cc}
    1 & a \\
    0 & b \\
  \end{array}
\right|=b\neq 0.
$$
Thus any $z\in\A$ is linear with respect to $x$.

Now we will describe the structure of module of $k$-forms $\A^k[\tau]$. We choose $(dx)^k$ as the basis for the right $\A$-module $\A^k[\tau]$, then any $k$-form $\omega$ can be written $\omega=(dx)^k\,u, \;u\in\A$. We have the following relations
$$
(dx)^k=\tau^k\,Q_k(x),\;\;\;d^kx=\tau^k\,P_k(x).
$$
In order to get a formula for the exterior differential of a $k$-form $\omega$ we need the polynomials $\Phi_1(x),\Phi_2(x),\ldots,\Phi_{N-1}(x)$ which can be defined by the recurrent relation
\begin{equation}
\Phi_{k+1}(x)=\text{Ad}_{\Delta x}(\Phi_{k}) + q^{k-1}\Phi_1(x),\quad k=1,2,\ldots,N-1,
\end{equation}
where $\Phi_1(x)=Q^{-1}_2(x)P_2(x)$. These polynomials satisfy the relations $d(dx)^k=(dx)^{k+1}\Phi_k(x)$ and given a $k$-form $\omega=(dx)^k\,u, \;u\in\A$ we find its exterior differential as
$$
d\omega=(dx)^{k+1}\,\bigg(q^k\,u^\prime_x + \Phi_k(x)\,u\bigg)=(dx)^{k+1}\,D^{(k)}u.
$$
The linear operator $D^{(k)}:\A\to\A, k=1,2,\ldots,N-1$ introduced in the previous formula has the form
\begin{equation}
D^{(k)}u=q^k\,u^\prime_x + \Phi_k(x)\,u,
\label{covariant derivative 2}
\end{equation}
and, as it was mentioned before, this operator resembles a covariant derivative of classical differential geometry. It is easy to see that the operator (\ref{covariant derivative 1}) is the particular case of (\ref{covariant derivative 2}), i.e. $D^{(1)}\equiv D$.
\section{Semi-Commutative Galois Extension Approach to Reduced Quantum Plane }
In this section we show that a reduced quantum plane can be considered as a semi-commutative Galois extension. We study a first order and higher order differential calculus of a semi-commutative Galois extension in the particular case of a reduced quantum plane.

\vskip.3cm
\noindent
Let $x,y$ be two variables which obey the commutation relation
\begin{equation}
x\,y=q\;y\,x,
\label{commutation relation on quantum plane}
\end{equation}
where $q\neq 0,1$ is a complex number. These two variables generate the algebra of polynomials over the complex numbers. This algebra is an associative algebra of polynomials over $\mathbb C$ and the identity element of this algebra will be denoted by $\1$. In noncommutative geometry and theoretical physics a polynomial of this algebra is interpreted as a function of {a quantum plane}
with two noncommuting coordinate functions $x,y$ and the algebra of polynomials is interpreted as the algebra of (polynomial) functions  of a quantum plane. If we fix an integer $N\geq 2$ and impose the additional condition
\begin{equation}
x^N=y^N=\1,
\label{x^N=y^N=1}
\end{equation}
then a quantum plane is referred to as {a reduced quantum plane} and this polynomial algebra will be denoted by $\reducedqp$.
\vskip.3cm
\noindent
Let us mention that from an algebraic point of view an algebra of functions on a reduced quantum plane may be identified with the generalized Clifford algebra ${\frak C}^N_2$ with two generators $x,y$. Indeed a generalized Clifford algebra
is an associative unital algebra generated by variables $x_1,x_2,\ldots,x_p$ obeying
the relations $x_ix_j=q^{\mbox{sg}(j-i)}x_jx_i,\; x_i^N=1$, where
$\mbox{sg}$ is the sign function.
\vskip.3cm
\noindent
It is well known that the generalized Clifford algebras have matrix representations, and, in the particular case of the algebra $\reducedqp$, the generators of this algebra $x,y$ can be identified with the square matrices of order $N$
\begin{equation}
x=\left(
    \begin{array}{cccccc}
      1 & 0 & 0 & \ldots & 0 & 0 \\
      0 & q^{-1} & 0 & \ldots & 0 & 0 \\
      0 & 0 & q^{-2} & \ldots & 0 & 0 \\
      \vdots & \vdots & \vdots & \ddots & \vdots & \vdots \\
      0 & 0 & 0 & \ldots & q^{-(N-2)} & 0 \\
      0 & 0 & 0 & \ldots & 0 & q^{-(N-1)} \\
    \end{array}
  \right),\;\;\;
y=\left(
    \begin{array}{cccccc}
      0 & 1 & 0 & \ldots & 0 & 0 \\
      0 & 0 & 1 & \ldots & 0 & 0 \\
      0 & 0 & 0 & \ldots & 0 & 0 \\
      \vdots & \vdots & \vdots & \ddots & \vdots & \vdots \\
      0 & 0 & 0 & \ldots & 0 & 1 \\
      1 & 0 & 0 & \ldots & 0 & 0 \\
    \end{array}
  \right),\label{matrix y}
\end{equation}
where $q$ is a primitive $N$th root of unity. As the matrices (\ref{matrix y}) generate the algebra $\mbox{Mat}_N({\mathbb C})$ of square matrices of order $N$ we can identify the algebra of functions on a reduced quantum plane with the algebra of matrices $\mbox{Mat}_N({\mathbb C})$.
\vskip.3cm
\noindent
The set of monomials $B=\{{\bf
1},y,x,x^2,yx,y^2,\ldots,y^kx^l,\ldots,y^{N-1}x^{N-1}\}$ can be
taken as the basis for the vector space of the algebra
$\reducedqp$. We can endow this vector
space with an ${\mathbb Z}_N$-graded structure if we assign degree zero to the identity element $\1$ and variable $x$ and we assign degree one to the variable $y$. As usual we define the degree of a product of two variables $x,y$ as the sum of degrees of factors. Then a polynomial
\begin{equation}
w=\sum_{l=0}^{N-1} \beta_{l}y^kx^l,\quad \beta_l\in {\mathbb C}, %
\label{homogeneous polynomial}
\end{equation}
will be a homogeneous polynomial with degree $k$.
Let us denote the degree of
a homogeneous polynomial $w$ by $|w|$ and the subspace of the
homogeneous polynomials of degree $k$ by $\reducedqpk$. It is
obvious that
\begin{equation}
\reducedqp=\reducedqpzero\oplus \reducedqpone
  \oplus\ldots\oplus {{\mathscr A}^{N-1}_q[x,y]}.%
  \label{graded structure}
\end{equation}
In particular a polynomial $r$ of degree zero can be written as follows
\begin{equation}
r=\sum_{l=0}^{N-1} \beta_lx^l,\qquad \beta_l\in {\mathbb
C},\;\;r\in\reducedqpzero.%
\label{x-polynomial}
\end{equation}
Obviously the subspace of elements of degree zero $\reducedqpzero$ is the subalgebra of $\reducedqp$ generated by the variable $x$.
Evidently the polynomial algebra $\reducedqp$ of polynomials of a reduced quantum plane can be considered as a semi-commutative Galois extension of the subalgebra $\reducedqpzero$ by means of the element $y$ which satisfies the relation $y^N=\1$. The commutation relation $x\,y=q\;y\,x$ gives us a semi-commutativity of this extension.

\vskip.3cm
\noindent
Now we can endow the polynomial algebra $\reducedqp$ with an $N$-differential $d$. Making use of Theorem 1 we define the $N$-differential by the following formula
\begin{equation}
dw=[y,w]_q=y\,w-q^{|w|}\,w\,y,
\end{equation}
where $q$ is a primitive $N$th root of unity and
$w\in \reducedqp$. Hence the algebra $\reducedqp$ equipped with the $N$-differential $d$ is a graded $q$-differential algebra.
\vskip.3cm
\noindent
In order to give a differential-geometric interpretation to the
graded $q$-differential algebra structure of $\reducedqp$ induced
by the $N$-differential $d_v$ we interpret the commutative
subalgebra $\reducedqpzero$ of the $x$-polynomials
(\ref{x-polynomial}) of $\reducedqp$ as an algebra of polynomial
functions on a one dimensional space with coordinate $x$. Since
$\reducedqpk$ for $k>0$ is a $\reducedqpzero$-bimodule we
interpret this $\reducedqpzero$-bimodule of the elements of
degree $k$ as a bimodule of differential forms of degree $k$ and
we shall call an element of this bimodule a differential $k$-form
on a one dimensional space with coordinate $x$. The
$N$-differential $d$ can be interpreted as an exterior differential.
\vskip.3cm
\noindent
It is easy to show that in one dimensional case we have a simple
situation when every bimodule $\reducedqpk, k>0$ of the
differential $k$-forms is a free right module over the commutative
algebra of functions $\reducedqpzero$. Indeed if we write a
differential $k$-form $w$ as follows
\begin{equation}
w=y^k\sum_{l=0}^{N-1} \beta_lx^l=y^k\,r,\quad%
r=\sum_{l=0}^{N-1}\beta_lx^l\in\reducedqpzero,
\end{equation}
and take into account that the polynomial $r=(y^k)^{-1}w=y^{N-k}w$
is uniquely determined then we can conclude that $\reducedqpk$ is
a free right module over $\reducedqpzero$ generated by $y^k$.
\vskip.3cm
\noindent
As it was mentioned before a bimodule structure of a free right module
over an algebra $\cal B$ generated freely by $p$ generators is
uniquely determined by the homomorphism from an algebra $\cal B$
to the algebra of $(p\times p)$-matrices over $\cal B$. In
the case of a reduced quantum plane every right module
$\reducedqpk$ is freely generated by one generator (for instant we
can take $y^k$ as a generator of this module). Thus its bimodule
structure induces an endomorphism of the algebra of functions
$\reducedqpzero$ and denoting this endomorphism in the case of the
generator $y^k$ by
$A_k:\reducedqpzero\to\reducedqpzero$ we get%
\begin{equation}
r\,y^k=y^k\,A_{k}(r), \;\;(\mbox{no summation over}\;\; k)
\label{bimodule structure}
\end{equation}%
for any function $r\in\reducedqpzero$. Making use of the commutation relations of variables $x,y$
we easily find that $A_{k}(x)=q^k\,x$. Since the
algebra of functions $\reducedqpzero$ may be viewed as a bimodule
over the same algebra we can consider the functions as degree zero
differential forms, and the corresponding endomorphism is the
identity mapping of $\reducedqp$, i.e. $A_0=I$, where
$I:\reducedqpzero\to \reducedqpzero$ is the identity mapping. Thus
the bimodule structures of the free right modules
$\reducedqpzero,\reducedqpone,\ldots,{{\mathscr A}^{N-1}_q[x,y]}$ of
differential forms induce the associated endomorphisms
$A_0,A_1,\ldots, A_{N-1}$ of the algebra $\reducedqpzero$. It is
easy to see that for any $k$ it holds $A_k=A^k_1$.
\vskip.3cm
\noindent
Let us start with the first order differential calculus $(\reducedqpzero,\reducedqpone,d)$ over the
algebra of functions $\reducedqpzero$ induced by the
$N$-differential $d$,
where $d:\reducedqpzero\to\reducedqpone$ and $\reducedqpone$ is
the bimodule over $\reducedqpzero$. For any $w\in\reducedqpzero$ we have
\begin{eqnarray}
dw=yw-wy=yw-yA_1(w)=y(w-A_1(w))=y\,\Delta_q(w),%
\end{eqnarray}
where $\Delta_q=I-A_1:\reducedqpzero\to \reducedqpzero$. It is
easy to verify that for any two functions $w,w'\in \reducedqpzero$ the
mapping $\Delta_q$ has the following properties
\begin{gather}
\Delta_q(ww')=\Delta_q(w)w'+A_1(w)\Delta_q(w'),%
                           \label{property1 of delta}\\%
\Delta_q(x^k)=(1-q)[k]_q\;x^k.%
    \label{propert2 of delta}
\end{gather}
\noindent
Particularly $dx=y\Delta_q(x)$, and this formula shows that
$dx$ can be taken as a generator for the free right module
$\reducedqpone$.
\vskip.3cm
\noindent
Since the bimodule $\reducedqpone$ of the first order differential
calculus $(\reducedqpzero,\reducedqpone,d)$ is a free right module we have
a coordinate first order differential calculus over the algebra $\reducedqpzero$, and in the case of a calculus of this
kind the differential induces the derivative
$\partial:\reducedqpzero\to\reducedqpzero$ which is defined by the
formula $dw=dx\,\,\partial w,\; \forall w\in \reducedqpzero$.
Using this definition we find that for any function $w$ it holds
\begin{equation}%
\partial w=(1-q)^{-1}x^{N-1}\Delta_q(w).%
\end{equation}%
From this formula and (\ref{property1 of delta}),(\ref{propert2 of
delta}) it follows that this derivative satisfies the twisted
Leibniz rule
\begin{equation}%
\partial(ww')=\partial(w)\cdot w'+A_1(w)\cdot\partial(w'),
\end{equation} and
\begin{equation}\partial x^k=[k]_q\;x^{k-1}.\end{equation}
\noindent
Let us study the structure of the higher order exterior calculus
on a reduced quantum plane or, by other words, the structure of
the bimodule $\reducedqpk$ of differential $k$-forms, when $k>1$.
In this case we have a choice for the generator of the free right
module. Indeed since the $k$th power of the exterior differential
$d$ is not equal to zero when $k<N$, i.e. $d^k\neq 0$ for
$k<N$, a differential $k$-form $w$ may be expressed either by
means of $(dx)^k$ or by means of $d^kx$. Straightforward
calculation shows that we have the following relation between
these generators
\begin{equation}
d^kx=\frac{[k]_q}{q^{\frac{k(k-1)}{2}}}(dx)^k\,x^{1-k}.
\end{equation}
We will use the generator $(dx)^k$ of the free right module
$\reducedqpk$ as a basis in our calculations with differential
$k$-forms. For any differential $k$-form $w\in\reducedqpk$ we
have $dw\in\reducedqpkplusone$. Let us express these two
differential forms in terms of the generators of the modules
$\reducedqpk$ and $\reducedqpkplusone$. We have $w=(dx)^k\,r,\,
dw=(dx)^{k+1}\,{\tilde r}$, where $r,{\tilde
r}\in\reducedqpzero$ are the functions. Making use of the
definition of the exterior differential $d$ we calculate the
relation between the functions $r,\tilde r$ which is
\begin{equation}
\tilde r=(\Delta_qx)^{-1}(q^{-k}r-q^kA_1(r)),
\end{equation}
where $A_1$ is the endomorphism of the algebra of functions
$\reducedqpzero$. This relation shows that the exterior
differential $d$ considered in the case of the differential
$k$-forms induces the mapping
$\Delta^{(k)}_q:\reducedqpzero\to\reducedqpzero$ of the algebra of
the function which is defined by the formula%
\begin{equation}%
dw=(dx)^{k+1}\Delta^{(k)}_q(r),
\end{equation}
where%
\begin{equation}
w=(dx)^k\,r.
\end{equation}
It is obvious that%
\begin{equation}
\Delta^{(k)}_q(r)=(\Delta_qx)^{-1}(q^{-k}r-q^kA_1(r)).
\end{equation}%
It is obvious that for $k=0$ the mapping $\Delta^{(0)}_q$
coincides with the derivative induced by the differential
$d$ in the first order calculus, i.e.%
\begin{equation}%
\Delta^{(0)}_q(r)=\partial r=(\Delta_q x)^{-1}(r-A_1(r)).
\end{equation}
The higher order mappings $\Delta_q^{(k)}$, which we do not have
in the case of a classical exterior calculus on a one dimensional
space, have the derivation type property
\begin{equation}
\Delta_q^{(k)}(r\,r')=\Delta_q^{(k)}(r)\,r'+q^k\,A_1(r)\,\Delta^{(0)}_x(r'),
\end{equation}
where $k=0,1,2,\ldots,N-1$. A higher order mapping
$\Delta_q^{(k)}$ can be expressed in terms of the derivative
$\partial$ as a differential operator on the algebra of functions
as follows
\begin{equation}
\Delta_q^{(k)}=q^k\,\partial\; +\frac{q^{-k}-q^k}{1-q}\,x^{-1}.
\label{deformed derivative}
\end{equation}
\noindent
Thus we see that exterior calculus on a one dimensional space with
coordinate $x$ satisfying $x^N=1$ generated by the exterior
differential $d$ satisfying $d^N=0$ has the differential forms
of higher order which are not presented in the case of a classical
exterior calculus with $d^2=0$. The formula for the exterior
differential of differential forms can be defined by means of contains not an
a derivative which satisfies the twisted Leibniz rule (\ref{deformed
derivative}).

\subsection*{Acknowledgment}
The authors is gratefully acknowledge the Estonian Science Foundation for financial support of this work under the Research Grant No. ETF9328. This research was also supported by institutional research funding IUT20-57
of the Estonian Ministry of Education and Research. The second author also acknowledges his gratitude to the Doctoral School in Mathematics and Statistics for financial support of his doctoral studies at the Institute of Mathematics, University of Tartu.

\end{document}